\documentclass[12pt]{article}

\font\tenmsb=msbm10
\font\sevenmsb=msbm7
\font\fivemsb=msbm5

\newfam\msbfam
\textfont\msbfam=\tenmsb
\scriptfont\msbfam=\sevenmsb
\scriptscriptfont\msbfam=\fivemsb
\def\Bbb#1{{\fam\msbfam #1}}

\makeatletter
\ifnum\@ptsize=0  \fi \ifnum\@ptsize=2
 \fi 

\newcommand\bC{{\Bbb C}}

\newcommand\bP{{\Bbb P}}

\newtheorem{theorem}{Th\' eor\`eme}[section]
\newtheorem{lemma}[theorem]{Lemme}
\newtheorem{corollary}[theorem]{Corollaire}
\newtheorem{proposition}[theorem]{Proposition}
\newtheorem{question}[theorem]{Question}
\newtheorem{re}[theorem]{Remarque}

\newtheorem{definition}[theorem]{D\'efinition}

\newtheorem{example}[theorem]{Exemple}

\newenvironment{remark}{\begin{re}\em}{\end{re}}

\begin{document}

\title {ISOTRIVIALIT\'E DE CERTAINES FAMILLES K\" AHL\'ERIENNES  DE VARI\'ET\'ES NON PROJECTIVES} 
\author{Fr\'ed\'eric Campana}

\maketitle

\section{Introduction}

\
Soit $f:X\to S$ est une application holomorphe surjective connexe avec $X$ K\" ahl\'erienne compacte et connexe, et soit $X_s$ une fibre lisse non projective de $f$. Nous montrons que $f$ est isotriviale si $X_s$ est soit hyperk\" al\'erienne irr\'eductible, soit un tore complexe g\'en\'eral, soit la vari\'et\'e de Kummer associ\'ee \`a un tel tore. La d\'emonstration repose sur le th\'eor\`eme de Torelli local pour les $2$-formes, et sur le crit\`ere de projectivit\'e de Kodaira: $X_s$ est projective si elle n'a pas de $2$-forme holomorphe non-nulle.

Nous commen\c cons par rappeler les notions tr\`es classiques utilis\'ees.

\subsection{Tores simples et vari\'et\'es hyperk\" ahl\'eriennes.} 

\

Soit $F$ une vari\'et\'e (lisse) K\"ahl\'erienne compacte et connexe de dimension complexe $n$.

On dit que $F$ est un tore complexe simple si $F$ n'a pas de sous-tore complexe non-trivial. On introduira en \ref{torir} ci-dessous la notion de tore irr\'eductible en poids $2$. Ces tores sont g\'en\'eraux dans leur famille de d\'eformation (voir \S\ref{torgenirr}).

La vari\'et\'e de Kummer $F'$ de $F$ est alors le quotient de l'\'eclat\'e $\hat{F}$ de $F$ en ses $2^{2n}$ points de torsion par l'involution de $\hat{F}$ qui rel\`eve l'involution $x\to (-x)$ de $F$. On dit alors que $F'$ est une vari\'et\'e de Kummer simple (resp. irr\'eductible en poids $2$) si le tore $F$ l'est. Voir [U 75], \S. 16,pp. 190-203 pour une \'etude d\'etaill\'ee de ces vari\'et\'es. En particulier, lorsque $n\geq 3$, les petites d\'eformations d'une vari\'et\'e de Kummer sont des vari\'et\'es de Kummer. (Par contre, lorsque $n=2$, $F'$ est une surface $K3$).

Rappelons que $F$ est  {\it hyperk\"ahl\'erienne}  si $F$ est simplement connexe de dimension complexe paire $n=2m$ et admet une $2$-forme holomorphe $u$ de rang maximum $m$ en chaque point. On dit que $F$ est {\it hyperk\"ahl\'erienne irr\'eductible} si, de plus,  $h^{2,0}(F)=1$. 

Toute vari\'et\'e hyperk\"ahl\'erienne se d\'ecompose en produit de vari\'et\'es hyperk\"ahl\'eriennes irr\'eductibles ([Bo74], [Be83] et [Y78]). Cette d\'ecomposition est unique ([Be83], th\'eor\`eme 2,p. 764). 

On dira que la vari\'et\'e hyperk\"ahl\'erienne $F$ est {\it sans facteur alg\'ebrique} si aucun des facteurs de cette d\'ecomposition n'est une vari\'et\'e projective.  On trouvera des \'etudes approfondies et des r\'ef\'erences r\'ecentes sur les vari\'et\'es hyperk\" ahl\'eriennes dans [H03] et [V02] .

\subsection{Fibrations.}

\

On appellera {\it fibration} toute application holomorphe surjective \`a fibres connexes $f:X\to S$ entre la vari\'et\'e analytique complexe compacte et k\"ahl\'erienne $X$ et l'espace analytique complexe $S$. On notera $X_s$ une fibre {\it g\'en\'erique lisse} de $f$, au-dessus de $s\in S$. 

On dit que $f$ est {\it isotriviale} (resp. {\it triviale}) si $X_s$ et $X_{s'}$ sont analytiquement isomorphes pour $s,s'$ g\'en\'eriques dans $S$ (resp. s'il existe une application bim\'eromorphe $t: S\times F\to X$ au-dessus de $S$ qui induit par restriction un isomorphisme de $F$ sur $X_s$ pour $s\in S$ g\'en\'erique). 

Si le groupe d'automorphismes analytiques de $X_s$ est discret, alors $f$ est isotriviale si et seulement s'il existe un changement de base $u: S'\to S$ g\'en\'eriquement fini tel que le morphisme $f':X':=(X\times_{S}S')\to S'$ d\'eduit de $f$ par $u$ soit trivial. 

\subsection{Le r\'esultat principal.} 

\

L'objectif du pr\'esent texte est de d\'emontrer entre autres le r\'esultat suivant, qui r\'epond \`a, et a \'et\'e motiv\'e par, une question pos\'ee par Anand Pillay, que nous remercions vivement.

\begin{theorem}\label{isot} Soit $f:X\to S$ une application holomorphe surjective \`a fibres connexes entre la vari\'et\'e K\"ahler compacte $X$ et l'espace analytique complexe $S$. Soit $X_s$ une fibre lisse de $f$.
Alors $f$ est isotriviale dans chacun des trois cas suivants:
 
 1. $X_s$ est une vari\'et\'e hyperk\"ahl\'erienne sans facteur alg\'ebrique (ceci est r\'ealis\'e, par exemple si la dimension alg\'ebrique $a(X_s)$ est nulle).
 
 2. $X_s$ est un tore complexe non projectif et irr\'eductible en poids $2$.
 
 3. $X_s$ est la vari\'et\'e de Kummer associ\'ee \`a un tore complexe non projectif et irr\'eductible en poids $2$.
\end{theorem}

\ 

Le cas 1. \'etait \'etabli dans [Ca89] par une toute autre approche lorsque $S$ est de dimension alg\'ebrique $0$.

La d\'emonstration sera donn\'ee dans les \S 2-5.

\

\begin{remark} Si l'on ne suppose pas que $X$ est compacte K\"ahler, la conclusion peut \^etre en d\'efaut dans chacun des trois cas: l'espace des twisteurs $Z$ associ\'e \`a une famille de Calabi de structures complexes sur $F$, une surface $K3$ Ricci-plate et non projective, est une vari\'et\'e complexe compacte de dimension complexe $3$ munie d'une submersion holomorphe non isotriviale $f:Z\to \bP_1(\bC)$ dont les fibres sont des surfaces $K3$, donc hyperk\" ahl\'eriennes irr\'eductibles. Mais $Z\notin \cal C$, car $Z$ est rationellement connexe, donc $h^{(2,0)}(Z)=0$. Cependant, $Z$ n'est pas projective, et donc pas K\"ahler, car la fibre $F$ ne l'est pas. Or, un  th\'eor\`eme de Kodaira ( utilis\'e crucialement ci-dessous) affirme que si $X$ est K\"ahler non projective, on a $h^{2,0}(X)\neq 0$.

On peut construire de la m\^eme fa\c con des exemples similaires avec $X_s$ un tore complexe irr\'eductible en poids $2$ non projectif de dimension $2$, ou la surface de Kummer associ\'ee \`a un tel tore (Vari\'et\'es de Blanchard. Voir [U75], Example 9.9, p. 105).

Il existe, par contre, de nombreuses fibrations non isotriviales $f:X\to S$ avec $X$ projective et \`a fibres lisses hyperk\" ahl\'eriennes irr\'eductibles.

Remarquons enfin que les vari\'et\'es $X_s$ apparaissant dans l'\'enonc\'e de \ref{isot} sont g\'en\'erales dans leur famille de d\'eformation (sauf si $n=2$ dans le cas 3).

\end{remark}

\subsection{Motivation: les vari\'et\'es simples.}

\

Rappelons que $X_s$ (compacte K\" ahler connexe) est dite {\it simple} si elle n'est pas recouverte par ses sous-espaces analytiques complexes compacts de dimension  $p$, pour $0<p<n$. La structure des vari\'et\'es compactes K\"ahler arbitraires peut-\^etre r\'eduite, au moyen de fibrations appropri\'ees, \`a celles des vari\'et\'es de Moishezon et des vari\'et\'es simples. Le nombre de fibrations n\'ecessaires serait ind\'ependant de la dimension si la question suivante avait une r\'eponse affirmative.

\begin{question} Soit $f:X\to S$ holomorphe surjective et connexe, avec $X$ compacte K\"ahler. Si $X_s$ est {\it simple}, alors $f$ est-elle isotriviale?
\end{question}

Les tores complexes simples non projectifs, ainsi que leurs vari\'et\'es de Kummer associ\'ees sont simples. Les vari\'et\'es hyperk\"ahl\'eriennes irr\'eductibles {\it g\'en\'erales} (dans leurs familles de d\'eformations) semblent devoir \^etre simples (exemples: les surfaces $K3$ ou ab\'eliennes et les d\'eformations g\'en\'erales de leurs produits sym\'etriques d\'esingularis\'es). Le th\'eor\`eme \ref{isot} est donc un test pour cette question.

Inversement:

\begin{question} Soit $X$ une vari\'et\'e K\"ahl\'erienne compacte simple, telle que $q(X')=0$, pour tout  $X'$ g\'en\'eriquement fini sur $X$. Alors $X$ est-elle de dimension paire $n=2m$, et a-t-elle une unique (\`a homoth\'etie complexe pr\`es) $2$-forme holomorphe, cette forme \'etant  de rang maximum $2m$ sur un ouvert de Zariski non vide $U$ de $X$? De plus, $X$ a-t'elle un mod\`ele minimal (\`a singularit\'es terminales) ayant un fibr\'e canonique de torsion?
\end{question}

Remarquons que l'on ne sait pas exclure l' existence de telles $X$ en dimension $3$, \`a moins d'admettre l'existence d'une th\'eorie de Mori K\" ahl\'erienne en cette dimension ([Pe01]).

\

Ces questions sont aussi pos\'ees ind\'ependamment par A. Pillay, dont les motivations viennent de la th\'eorie des Mod\`eles (voir [P-S02]). De ce point de vue, il s'agit de savoir si la classe des vari\'et\'es K\"ahler compactes est {\it non-multidimensionnelle}, c'est-\`a-dire si lorsque $f:X\to Y$ est une fibration, avec $X$ compacte K\"ahler, et si $a,b\in Y$ sont g\'en\'eriques, il existe un sous-espace analytique compact irr\'eductible non-trivial de $X_a\times X_b$ se projettant surjectivement sur chacun des facteurs. Le cas g\'en\'eral peut \^etre r\'eduit au cas o\`u $X_a$ est {\it simple} pour $a\in Y$ g\'en\'eral.

\

\begin{remark} Le pr\'esent texte est une version consid\'erablement simplifi\'ee de [C04], o\`u le th\'eor\`eme \ref{isot} \'etait \'etabli dans le cas hyperk\" ahl\'erien par une m\'ethode compliqu\'ee et peu naturelle, faute de disposer de la proposition \ref{restr}, pos\'ee sous forme de question. La d\'emonstration donn\'ee ici est due \`a C. Voisin. L'\'enonc\'e \ref{restr} permet, en particulier, de traiter les tores complexes g\'en\'eraux. Nous la remercions vivement. \end{remark}

\section{Restriction des $2$-formes holomorphes.}\label {hodge}

\

\begin{proposition} \label{restr} ([Voi04]) Soit $f:X\to S$ une fibration, avec $X$ une vari\'et\'e k\" ahl\'erienne compacte connexe. Soit $X_s$ une fibre lisse g\'en\'erique de $f$. Si l'application de restriction $r:H^0(X,\Omega_X^2)\to H^0(X_s,\Omega_{X_s}^2)$ est nulle, alors $X_s$ est projective.
\end{proposition}

Cette proposition est une version relative du th\'eor\`eme de Kodaira [Kod54], d'apr\`es lequel $X_s$ est projective si $ H^0(X_s,\Omega_{X_s}^2)=0$. Elle pr\'esente un int\'er\^et ind\'ependant, et de nombreuses applications potentielles.

\

{\bf D\'emonstration:} L'application de restriction $r:H^2(X,\Bbb C)\to H^2(X_s,\Bbb C)$ est d\'efinie sur $\Bbb Q$ et pr\'eserve les types (ie: est compatible avec la d\'ecomposition de Hodge). Supposons que la restriction \`a $X_s$ de toute $2$-forme holomorphe $u$ sur $X$ soit nulle. L'image de $H^2(X,\Bbb R)$ est donc contenue dans $H^{1,1}(X_s,\Bbb R)$ et coincide avec celle de $H^{1,1}(X,\Bbb R)$. Mais $H^{1,1}(X,\Bbb R)$ contient le c\^one ouvert non vide des classes de K\" ahler. Par densit\'e, il existe donc $v\in H^2(X,\Bbb Q)$ ayant m\^eme image par $r$ qu'une classe de k\" ahler $w$ sur $X$. La restriction de $w$ \`a $X_s$ est donc une classe de K\" ahler rationnelle sur $X_s$. Le th\'eor\`eme de Kodaira montre que $X_s$ est projective.

\

Sous des hypoth\`eses plus restrictives, on obtient la non-nullit\'e de {\it toutes} les restrictions:

\begin{proposition}\label{restr'} Soit $f:X\to C$ une fibration, avec $X$ compacte K\"ahler, et $C$ une courbe complexe projective. Si $X_s$ est une fibre lisse de $f$, et si $b_1(X_s)=0$, toute $2$-forme holomorphe non-nulle $u$ sur $X$ a une restriction $u_s$ \`a $X_s$ qui est non-nulle.
\end{proposition}

{\bf D\'emonstration:} Dans des coordonn\'ees locales $(t,z)$ sur $X$, o\`u $t$ est une coordonn\'ee sur $C$, et $z=(z_1,\dots,z_{n-1})$ des coordonn\'ees sur $X_s$, on peut \'ecrire sur $X_s$: $0\neq u=dt\wedge w+u_s$, o\`u $w$ est une $1$-forme intrins\`equement d\'efinie sur $X_s$. Donc $w=0$, et $u_s\neq 0$.

\section{Crit\`ere d'isotrivialit\'e.}\label{torel}

\

\begin{definition} On dit que $X$, compacte K\"ahler, satisfait le th\'eor\`eme de Torelli local pour les $2$-formes holomorphes si les p\'eriodes des $2$-formes holomorphes sur $X$ d\'eterminent localement dans son espace de Kuranishi la structure analytique de $X$. 
\end{definition}

\

\begin{example}\label{torloc} Les vari\'et\'es hyperk\" ahl\'eriennes irr\'eductibles satisfont le th\'eor\`eme de Torelli local pour les $2$-formes holomorphes, par [Be83], ainsi que les tores complexes et leurs vari\'et\'es de Kummer.
\end{example} 

\

\begin{definition} \label{torir} Soit $X$ une vari\'et\'e K\" ahl\'erienne compacte. On dit que $X$ est {\it irr\'eductible en poids $2$} si toute sous-structure de Hodge rationnelle $V:=V_{\Bbb Q}\otimes \Bbb C$ de la structure de Hodge sur $H^2(X,\Bbb C)$ contient $H^{(2,0)}(X)$, si $0\neq V^{(2,0)}:=V\cap H^{(2,0)}(X)$ .

\end{definition}

\

\begin{example} \label{irredtor} 

\

0. Si $X$ est irr\'eductible en poids $2$, elle n'admet pas de fibration holomorphe non-triviale $f:X\to S$ sur une vari\'et\'e k\" ahl\'erienne $S$ telle que $h^{(2,0)}(S)\neq 0$. Un tore irr\'eductible en poids $2$ est donc simple, s'il n'admet pas de tore quotient de dimension $1$ (donc s'il est de dimension alg\'ebrique $0$).

1. Si $H^{(2,0)}(X)=1$, et en particulier si $X$ est hyperk\" ahl\'erienne irr\'eductible, alors $X$ est irr\'eductible en poids $2$. 

2. Un tore complexe $X$ complexe de dimension $2$, ainsi que sa vari\'et\'e de Kummer $X'$ associ\'ee sont irr\'eductibles en poids $2$.

3. En dimension $3$ ou plus, un tore complexe g\'en\'eral, ainsi que sa vari\'et\'e de Kummer associ\'ee sont irr\'eductibles en poids $2$. Ce r\'esultat sera \'etabli s\'epar\'ement dans le \S\ref{torgenirr}. Il serait int\'eressant de savoir si les tores simples sont irr\'eductibles en poids $2$.

\end{example}

\begin{proposition}\label{critisot} Soit $f:X\to S$ une fibration, avec $X$ k\" ahler, telle qu'une fibre lisse $X_s$ de $f$ soit irr\'eductible en poids $2$ et satisfasse le th\'eor\`eme de Torelli local pour les $2$-formes holomorphes. 

S'il existe, de plus, une $2$-forme holomorphe $u$ sur $X$ dont la restriction $u_s$ \`a $X_s$ est non-nulle, alors $f$ est isotriviale.

\end{proposition}

{\bf D\'emonstration:}
Par le th\'eor\`eme des cycles invariants de Deligne ([D71] et [D74]; voir [Voi02], \S16 pour une exposition particuli\`erement lisible), valable aussi dans le cas k\" ahler, $j^*(H^2(X,\Bbb C))$ est une sous-structure de Hodge rationnelle de $H^2(X_s,\Bbb C)$, si $j:X_s\to X$ est l'inclusion naturelle. Cette structure contient la classe de $u_s$, et coincide donc avec $H^2(X_s,\Bbb C)$, puisque $X_s$ est suppos\'ee irr\'eductible en poids $2$. Donc la restriction de $H^{(2,0)}(X)$ \`a $X_s$ coincide avec $H^{(2,0)}(X_s)$. Les formes $u$ sont $d$-ferm\'ees. Donc leurs p\'eriodes sur $X_t$ sont ind\'ependantes de $t\in S^*$ (dans une trivialisation diff\'erentiable locale de $f$), et \'egales \`a celles de $X_s$. Puisque $X_t$ satisfait le th\'eor\`eme de Torelli local, on en d\'eduit que $f$ est isotriviale sur $S^*$. On conclut, par exemple, par la compacit\'e de $Chow(X)$ ([L75]).

\

\begin{remark}\label{ssdelign} La d\'emonstration de la proposition \ref{critisot} ne fait pas usage du th\'eor\`eme des cycles  invariants de Deligne  si $h^{(2,0)}(X_s)=1$, et en particulier si $X_s$ est hyperk\" ahler irr\'eductible. \end{remark}

\

De la proposition \ref{critisot} et des exemples \ref{irredtor} et \ref{torloc}, on d\'eduit le:

\begin{corollary} \label{torelli}Soit $f:X\to S$ une fibration avec $X$ k\" ahler, et $X_s$ une fibre lisse de $f$. On suppose qu'il existe sur $X$ une $2$-forme holomorphe dont la restriction \`a $X_s$ est non-nulle. Alors $f$ est isotriviale dans chacun des trois cas suivants:

1. $X_s$ est une vari\'et\'e hyperk\"ahl\'erienne irr\'eductible.
 
 2. $X_s$ est un tore complexe irr\'eductible en poids $2$.
 
 3. $X_s$ est la vari\'et\'e de Kummer associ\'ee \`a un tore complexe irr\'eductible en poids $2$.

\end{corollary}

\

Le th\'eor\`eme \ref{isot} en r\'esulte imm\'ediatement, gr\^ace \`a \ref{restr}, sauf lorsque $X_s$ est hyperk\"ahler sans facteur alg\'ebrique, mais n'est pas irr\'eductible. Ce cas sera trait\'e dans le \S.\ref{concl} ci-dessous.

\section{Tores irr\'eductibles en poids $2$.}\label{torgenirr}

Rappelons que si $S$ est un espace analytique complexe, un point de $S$ est dit g\'en\'eral s'il appartient \`a  une intersection d\'enombrable d'ouverts de Zariski denses de $S$. Une vari\'et\'e complexe compacte et connexe est dite g\'en\'erale dans sa famille de d\'eformation si elle l'est dans sa famille de Kuranishi.

La proposition qui suit est certainement bien connue, nous en donnons une preuve, faute d'avoir trouv\'e une r\'ef\'erence ad\'equate (bien qu'un \'enonc\'e analogue pour les vari\'et\'es ab\'eliennes polaris\'ees se trouve dans [B-L,17.4.1-3]).

\begin{proposition} Les tores complexes irr\'eductibles en poids $2$ de dimension $m\geq 3$ sont g\'en\'eraux dans leur famille de d\'eformation. Plus pr\'ecis\'ement: le tore g\'en\'eral $T$ de dimension $m\geq 3$ n'a pas de sous-structure de Hodge rationnelle non-triviale de $H^2(T,\Bbb C)$.

\end{proposition}

{\bf D\'emonstration:} (Voir aussi la remarque 4.2 ci-dessous)
On consid\'erera un tore $T$ comme la donn\'ee d'un $\Bbb Q$-espace vectoriel $V$ de dimension $2m$, et d'une structure complexe $J\in End_{\Bbb R} (V_{\Bbb R})$. L'ensemble $C$ des structures complexes $J$ sera identifi\'e \`a la Grassmannienne des $m$-plans complexes $H=H_J^{(0,1)}\subset V_{\Bbb C}=V_{\Bbb R}\oplus i.V_{\Bbb R}$ tels que: $H\cap V_{\Bbb R}=H\cap i.V_{\Bbb R}=\lbrace 0\rbrace$.

Par d\'enombrabilit\'e de $\Bbb Q$, il suffit de d\'emontrer que si $\lbrace 0 \rbrace\neq W\subset \wedge_{\Bbb Q}^2V$ est un sous-espace vectoriel (d\'efini sur $\Bbb Q$, donc), et si $W_{\Bbb C}$ est une sous-structure de Hodge de $(\wedge^2V)_{\Bbb C}$ pour tout $J\in C$, alors $W=\wedge_{\Bbb Q}^2V$. C'est ce que nous allons maintenant d\'emontrer.

Si $w\in W$, on notera $r(w)$ le rang de $w$, c'est-\`a-dire le plus grand des entiers $s\geq 0$ tels que $w^{\wedge s}\neq 0$. On a, bien s\^ur: $0\leq r(w)\leq m$.

On distingue $2$ cas:

$\bullet$ Premier cas: il existe $w\in W$ tel que $r(w)=m$.

Dans une $\Bbb Q$-base ad\'equate $\lbrace dx_j,dy_j,j=1,\dots,m\rbrace$ de $V$, on a donc: $w=\sum_{j=1}^{j=m}dx_j\wedge dy_j$. Soit $J_0\in C$ d\'efinie par: $J_0.dx_j:=dy_j,\forall j$. On a donc: $2w=i.\sum_{j=1}^{j=m}dz_j\wedge d\bar{z}_j$, avec $dz_j:=dx_j+i.dy_j,\forall j$.

L'espace $TC_{J_0}$ tangent \`a $C$ en $J_0$ est canoniquement identifi\'e \`a $Hom_{\Bbb C}(H,\bar{H})$, si $H:=H_{J_0}^{(0,1)}$. On a donc, si $t\in TC_{J_0}$, si $s\in \Bbb C$, et si $J_s:\equiv J_0+s.t$ (en notant $\equiv$ l'\'egalit\'e modulo $O(\vert s\vert ^2$):

$d\bar{z}_j=d\bar{z}'_j-s.\sum_{j=1}^{j=m}t_{j,h}dz_h$, ceci pour tout $j$, si les $\lbrace d\bar{z}'_j,\forall j\rbrace$ forment une $\Bbb C$-base de $H_{J_s}^{(0,1)}$. Donc, par un calcul imm\'ediat: $$2w\equiv i.\sum_{j=1}^{j=m}dz'_j\wedge d\bar{z}'_j-i.s.\sum_{1\leq j<l\leq m}(t_{j,l}-t_{l,j}).dz_j\wedge dz_l+i.\bar s.\sum_{1\leq j<l\leq m}(\bar t_{j,l}-\bar t_{l,j}).d\bar z_j\wedge d\bar z_l.$$

Puisque, par hypoth\`ese, $W_{\Bbb C}$ est une sous-structure de Hodge de $(\wedge ^2V)_{\Bbb C}$ pour tout $J$ de $C$, on a: $\sum_{1\leq j<l\leq m}(t_{j,l}-t_{l,j}).dz_j\wedge dz_l\in W_{\Bbb C}$, pour tout choix de $t\in TC_{J_0}$. 

Donc: $$(*) Re(dz_j\wedge dz_l)=dx_j\wedge dx_l-dy_j\wedge dy_l\in W$$ et $$(**) Im(dz_j\wedge dz_l)=dx_j\wedge dy_l-dx_l\wedge dy_j\in W,$$ ceci pour tous $1\leq j<l\leq m$.

On en d\'eduit que $w':=w+(dx_k\wedge dy_l-dx_l\wedge dy_k)=\sum_{j=1}^{j=m}dx'_j\wedge dy'_j\in W$, avec: $dx'_j=dx_j,\forall j$, et $dy'_j=dy_j$ si $j\neq k,l$, et: $dy'_k:=dy_k+dy_l$, et enfin: $dy'_l=dy_l-dy_k$.

On donc, appliquant les relations $(*)$ et $(**)$ pr\'ec\'edentes \`a $w'$, et en tenant compte de $(*),(**)$:
$$(***)dx_j\wedge dy_k\in W, \forall j\neq k$$.

$$(****)dy_j\wedge dy_k\in W,\forall j\neq k$$.

Enfin (puisque $dx'_k\wedge dx'_l-dy'_k\wedge dy'_l\in W$ et $dx_k\wedge dx_l-dy_k\wedge dy_l\in W$):

$$(+)dx_k\wedge dy_k+dx_l\wedge dy_l\in W, \forall k\neq l.$$

D'o\`u l'on d\'eduit, puisque $m\geq 3$:

$$(++)dx_j\wedge dy_j\in W, \forall j$$.

On obtient les derni\`eres \'egalit\'es cherch\'ees: $$(+++)dx_j\wedge dx_k\in W, \forall j\neq k$$ en appliquant $(**)$, comme ci-dessus, \`a $w":=w+(dx_k\wedge dx_l-dy_k\wedge dy_l)$, pour tous $k,l$.

On a donc bien: $W=\wedge^2V$ dans ce cas.

$\bullet$ Second cas: $r(w)\leq (m-1)$ pour tout $w\in W$. On va montrer que ce cas ne peut pas se produire. Si $w\in W$ est de rang maximum $(r-1)\leq (m-1)$, on a, comme ci-dessus: $w=\sum_{j=1}^{j=r}dx_j\wedge dy_j$. 

Soit $J\in C$ d\'efinie par: $J.dx_j:=dy_j,\forall j=2,3,\dots,(r-1)$, tandis que $J.dx_1:=dx_{r}$, et $J.dy_1:=dy_{r}$. 

On d\'efinit ensuite: $dz_1:=dx_1+idx_r$, $dz_r:=dy_1+idy_r$, et $dz_j=dx_j+idy_j$ si $j=2,3,\dots, (r-1)$. 

On a alors: $w=Re(dz_1)\wedge Re(dz_r)+\sum_{j=2}^{j=(r-1)}Re(dz_j)\wedge Im(dz_j)$, et donc: $w^{(2,0)}=dz_1\wedge dz_r$. 

D'o\`u: $w':=Re(dz_1\wedge dz_r)=dx_1\wedge dy_1-dx_r\wedge dy_r\in W$. Donc: $w+w'\in W$. Or $r(w+w')=r$. Contradiction avec notre hypoth\`ese $r(w")\leq (r-1),\forall w"\in W$. Ce cas ne se produit donc pas.

\begin{remark} (due \`a C. Voisin) On peut donner des preuves plus courtes et conceptuelles de 4.1 ci-dessus en utilisant, soit la th\'eorie des variations infinit\'esimales de strctures de Hodge de Griffiths, soit l'action de la monodromie de la famille des tores complexes param\'etr\'ee par $C$.
\end{remark}

\

\section{Le cas hyperk\" ahler sans facteur alg\'ebrique.} \label{concl}

\

On va r\'eduire le cas $X_s$ hyperk\" ahler  du th\'eor\`eme \ref{isot} au cas hyperk\" ahl\' erien irr\'eductible gr\^ace au lemme suivant:

 \begin{lemma}\label{decrel} Soit $f:X\to S$ une fibration telle que $X_s$ soit hyperk\" ahl\'erienne irr\'eductible. Soit $s\in S^*$ g\'en\'erique. Alors $X_s$ se d\'ecompose de mani\`ere unique en un produit de facteurs hyperk\"ahl\'eriens irr\'eductibles $X_{j,s}; j=1,...,k$, dont le nombre $k$ de facteurs est constant (\'egal \`a $h^{(0,2)}(X_s)$). On note $r_{j,s}:X_s\to X_{j,s}$ la projection correspondante.

Il existe \footnote {apr\`es un changement de base fini $v:S'\to S$, notationnellement ignor\'e ici, pour simplifier.} des fibrations $r_j:X\to X_j$ et $f_j:X_j\to S, j=1,...,k$, telles que:
 
 1. $f=r_j\circ f_j, \forall j$.
 
 2. $X_{j,s}\cong (X_j)_s, \forall j$, et pour $s\in S$ g\'en\'erique.
 
 3. $(r_j)_{\vert X_s}:X_s\to X_{j,s}=r_{j,s}, \forall j$, et pour $s\in S$ g\'en\'erique.

 (Autrement dit: $X\cong X_1\times_{S}\times...\times _S\times X_k$ est produit fibr\'e au-dessus de $S$ de $k$ fibrations \`a fibres hyperk\" ahl\'eriennes irr\'eductibles, et cette d\'ecomposition induit sur $X_s$ la d\'ecomposition en les produits $X_{j,s}$).
 \end{lemma}
 
{ \bf D\'emonstration.} Reprenant les notations pr\'ec\'edentes, chacun des facteurs $X_{j,s}$ d\'etermine une projection $r_{j,s}:X_s\to X_{j,s}$ sur ce facteur. 

On va maintenant utiliser l'espace des cycles $\cal C$$(X)$ d'un espace analytique complexe $X$ (voir [B75]). Chacune des $r_{j,s}$ d\'etermine donc une unique composante connexe $C_{j,s}\cong X_{j,s}$ de $\cal C$$(X_s)$, celle qui consiste en la famille des fibres de cette projection, affect\'ees de  la multiplicit\'e $1$. La projection $r_{j,s}:X_s\to X_{j,s}$ est alors obtenue naturellement en consid\'erant le graphe $Z_{j,s}\subset C_{j,s}\times X_s$: ce graphe est isomorphe \`a $X_s$ par la (restriction \`a $Z_{j,s}$ de la) seconde projection $p_{j,s}$, et si $q_{j,s}$ est la premi\`ere projection, la projection de $X_s$ sur $X_{j,s}$ est simplement $q_{j,s}\circ (p_{j,s})^{-1}$. 

On va simplement relativiser cette construction au-dessus de $S$.

  Le th\'eor\`eme de Baire appliqu\'e \`a l'espace des cycles relatifs $f_*:\cal C$$(X/S)\to S$,  constitu\'e du sous-ensemble analytique ferm\'e $\cal C$$(X/S)$ de $\cal C$$(X)$, r\'eunion de tous les cycles $Z$ de $X$ dont le support est contenu dans une fibre $X_s$ de $f$, avec $s:=f_*(Z)$, montre qu'il existe, pour chaque $j=1,...,k$, une unique composante irr\'eductible $C_j$ de $\cal C$$(X/S)$ telle que $C_{j,s}$ soit la fibre de $C_j$ au-dessus de $s\in S$ g\'en\'erique \footnote {On remplace, pour chaque $j$, $S$ par la factorisation de Stein $S_j$ de $C_j$ au-dessus de $S$. On prend pour nouvelle base $S'$ une composante irr\'eductible principale du produit fibr\'e des $S_j$ au-dessus de $S$.}. Le graphe $Z_j\subset C_j\times _S X$ de $C_j$ est alors naturellement bim\'eromorphe \`a $X$ au-dessus de $S$, et lui est isomorphe au-dessus d'un ouvert de Zariski dense ad\'equat (par la seconde projection $p_j$). On d\'efinit alors $r_j:=q_j\circ (p_j)^{-1}:X\to C_j:=X_j$, $q_j:Z_j\to C_j$ \'etant simplement la seconde projection. De m\^eme, $f_j:C_j\to S$ est simplement la restriction de $f_*$. Les propri\'et\'es restantes se v\'erifient par restriction \`a $X_s$, $s\in S$ g\'en\'erique.
  
  \

 On obtient donc le cas $X_s$ hyperk\"ahler du th\'eor\`eme \ref{isot} lorsque $X$ est compacte k\" ahler en appliquant le lemme suivant \`a chacune des fibrations  $f_j:X_j\to S$, puisque par hypoth\`ese, leurs fibres g\'en\'erales sont hyperk\"ahl\'eriennes irr\'eductibles et non alg\'ebriques.
 
 \
 

\section{Bibliographie.}

\

[B75] D. Barlet. Espace analytique r\'eduit des cycles analytiques compacts d'un espace analytique complexe de dimension finie. LNM 482(1975), 1-158

[Be83] A. Beauville. Vari\'et\'es K\" ahl\'eriennes dont la premi\`ere classe de Chern est nulle.  J. Differential Geom. 18 (1983),755--782 .

[B-L 02] C.Birkenhake-H. Lange. Complex Abelian Varieties. Springer 1980 (Seconde \'edition 2002).

[Bo74] F. Bogomolov. The decomposition of K\" ahler manifolds with a trivial canonical class.  Mat. Sb. 93(135) (1974), 573--575.

[C89] F. Campana. Geometric algebraicity of moduli spaces of compact K\"  ahler symplectic manifolds. J. f. die reine u. angew. Math. 397 (1989), 202--207. 

[C04] F.Campana. Crit\`ere d'isotrivialit\'e pour les familles de vari\'et\'es hyperk\" ahl\'eriennes sans facteur alg\'ebrique. math.AG/0408148.

[D71] P. Deligne. Th\'eorie de Hodge II. Publ. IHES. 40 (1971), 5-57.

[D71] P. Deligne. Th\'eorie de Hodge III. Publ. IHES. 44 (1974), 5-77.

[H03] D. Huybrechts. Finiteness results for compact hyperk\" ahler manifolds. J. f. die reine u. angew. Math. 558 (2003), 15--22.

[Kod54] K. Kodaira. On K\" ahler varieties of restricted type. Ann. Math. 60 (1954), 28-48.

[L75] D. Lieberman. Compactness of the Chow scheme: applications to automorphisms and deformations of K\" ahler manifolds. LNM 670(1978), 140-186.

[Pe01] T. Peternell. Towards a Mori theory on compact K\" ahler threefolds. III. Bull. Soc. Math. France 129 (2001), 339--356.

[P-S02] A. Pillay- T. Scanlon. Compact complex manifolds with the DOP and other properties. Journal of Symbolic Logic, 67(2002), 737-743.

[U75] K. Ueno. Classification theory of algebraic varieties and compact complex manifolds. LNM 439 (1975).

[Ve02] M. Verbitsky. Hyperk\" ahler manifolds with torsion, supersymmetry and Hodge theory. Asian J. Math. 6 (2002), no. 4, 679--712.

[Voi04] C.Voisin. Th\'eorie de Hodge et g\'eom\'etrie alg\'ebrique complexe. Collections Cours sp\'ecialis\'es 10. SMF (2002).

[Voi04] C.Voisin. message electronique. 11/12/04.

[Y78] ST. Yau. On the Ricci curvature of a compact K\" ahler manifold and the complex Monge-Amp\`ere equation. I. Comm. Pure Appl. Math. 31 (1978), 339--411.

Ê

\section{Adresse}

\

F.Campana

Universit\'e Nancy 1.

D\'epartement de Math\'ematiques.

BP 239

F. 54506. Vandoeuvre-les-nancy. C\'edex.

campana@iecn.u-nancy.fr

\end{document}